\documentclass[
11pt]{amsart}

\usepackage{amsthm, amsfonts, amssymb}
\usepackage[pagebackref,colorlinks]{hyperref}
\usepackage[all]{xy}

\theoremstyle{definition}
\newtheorem{ntn}{Notation}[section]
\newtheorem{dfn}[ntn]{Definition}
\theoremstyle{plain}
\newtheorem{lem}[ntn]{Lemma}
\newtheorem{sublem}[ntn]{Sublemma}
\newtheorem{prp}[ntn]{Proposition}
\newtheorem{thm}[ntn]{Theorem}

\theoremstyle{remark}

\newtheorem{rem}[ntn]{Remark}
\newtheorem{exa}[ntn]{Example}

\def\floor[#1]{\lfloor #1 \rfloor }

\newcommand{\z}{\mathbb{Z}}

\newcommand{\q}{\mathbb{Q}}
\newcommand{\R}{\mathbb{R}}

\newcommand{\lan}{\langle}
\newcommand{\ran}{\rangle}

\newcommand{\GL}{\mathit{{\rm GL}}}

\newcommand{\SL}{\mathit{{\rm SL}}}

\newcommand{\SK}{\mathit{{\rm SK}}}

\newcommand{\ppp}{\mathfrak{p}}

\newcommand{\PP}{\mathcal{P}}

\renewcommand{\H}{\tilde{H}}

\newcommand{\inc}{{\rm inc}}

\newcommand{\tors}{{{\rm Tor}_1^{\z}}}

\newcommand{\arr}{\rightarrow}
\newcommand{\larr}{\longrightarrow}

\newcommand{\mt}{\mapsto}

\newcommand{\fff}{{F^\times}}

\newcommand{\rr}{{R^{\times}}}

\renewcommand{\char}{{\rm char}}

\renewcommand{\ker}{{\rm ker}}
\newcommand{\coker}{{\rm coker}}
\newcommand{\im}{{\rm im}}
\newcommand{\ind}{{\rm ind}}

\newtheoremstyle{athm}
  {}
  {}
  {\itshape}
  {}
  {\scshape}
  {}
  {.5em}
  {\thmnote{#3}}
\theoremstyle{athm}

\begin{document}

\title{A note on third homology of $\GL_2$}
\author{Behrooz Mirzaii}

\begin{abstract}
In this paper the third homology group of the linear group
$\GL_2(R)$ with integral coefficients is investigated,
where $R$ is a commutative ring with many units.
\end{abstract}

\maketitle

\section*{Introduction}

In this article the third homology group of $\GL_2(R)$ with
integral coefficients, i.e. $H_3(\GL_2(R), \z)$, is studied,
where $R$ is a commutative ring with many units, e.g. a semilocal ring with
infinite residue fields. Interest in the study of this group
comes from its relation to the Bloch group $B(R)$
\cite{suslin1991}, \cite{mirzaii-2009},
its connection to the scissor congruence problem in three
dimensional hyperbolic geometry \cite{dupont-sah1982}, etc.

A theorem of Bloch and Wigner claims the existence of the
exact sequence
\[
0 \arr \q/\z \arr H_3(\SL_2(F), \z) \arr \ppp(F)
\overset{\lambda}{\arr}
\bigwedge\!{}_\z^2 \fff \arr K_2(F) \arr 0,
\]
where $F$ is an algebraically closed field of $\char(F)=0$,
$\ppp(F)$ is the pre-Bloch group of $F$ and
$\lambda:[a] \mt a \wedge (1-a)$. We should mention that
$H_3(\SL_2(F), \z)$ is isomorph to $K_3(F)^\ind$,
where $K_3(F)^\ind$ is the indecomposable part of the third $K$-group
$K_3(F)$.

This theorem was generalized by Suslin \cite{suslin1991}to all infinite fields.
Suslin's methods are difficult and it is not clear how can one
generalize them to other nice rings,
for example to semilocal rings with infinite residue fields.
Recently, we were able to establish a version of Bloch-Wigner
exact sequence over rings with many units \cite{mirzaii-2009}.

The current paper grew out of our desire to understand the Bloch-Wigner
exact sequence, its proof and its connection to $K$-theory.
Here we try to shed some light to the role of the group
$H_3(\GL_2(R), \z)$ in this direction.
Let $R$ be a commutative ring with many units. Our main theorem in this paper
states the existence of the exact sequence
\begin{gather*}
0 \! \arr \! H_3(\GL_2(R),\z)/H_3(\GL_1(R),\z)\! \arr \! \PP(R)
\! \overset{\nu}{\arr}
\!\!\!\!\!\!\!\!\!\!\underset{(a_1,a_2) \neq (1,1)} {\coprod_{a_i \in \rr}}
\!\!\!\!\!\!\!\! \z.[a_1,a_2]
\! \arr \! K_2(R)\! \arr \! 0,
\end{gather*}
where $\PP(R)$ is a group closely related to the pre-Bloch group $\ppp(R)$.
We construct easily a natural map from  $\PP(R)$ to $\ppp(R)$
and using the above exact sequence show that under this map
$H_3(\GL_2(R), \z)/H_3(\GL_1(R),\z)$ maps
into the Bloch group $B(R):=\ker(\lambda)$. At the end we study a group
that is closely related to the Bloch group $B(R)$.

\subsection*{Notation}
In this paper by $H_i(G)$ we mean  the  homology of group $G$
with integral coefficients, namely $H_i(G, \z)$.
By $\GL_n$ (resp. $\SL_n$) we mean the general (resp. special)
linear group $\GL_n(R)$ (resp. $\SL_n(R)$).
If $A \arr A'$ is a homomorphism of abelian groups, by $A'/A$ we
mean $\coker(A \arr A')$. We denote an element of $A'/A$ represented
by $a' \in A'$ again by $a'$.

\section{Rings with many units}
In this article we always assume that rings are commutative and contain
the unit element $1$.
\begin{dfn}
We say that a commutative ring $R$ is a {\it ring with many units} if
for any $n \ge 2$ and for any finite number of surjective linear forms
$f_i: R^n \arr R$, there exists a $v \in  R^n$ such that, for all $i$,
$f_i(v) \in \rr$.
\end{dfn}

The study of rings with many units is originated by  W. van der
Kallen in \cite{vdkallen1977}, where he showed that $K_2$ of such
rings behaves very much like $K_2$ of fields.

For a ring $R$, the $n$-th Milnor $K$-group $K_n^M(R)$ is
defined as an abelian group
generated by symbols $\{a_1, \dots, a_n\}$, $a_i \in \rr$,
subject to multilinearity and the relation $\{a_1, \dots, a_n \}=0$
if there exits $i, j$, $i \neq j$, such that $a_i+a_j=0$ or  $1$.

Van der Kallen proved \cite{vdkallen1977} that when $R$ is a
ring with many units, then
\begin{equation}\label{k2}
K_2(R)\simeq K_2^M(R) \simeq \rr \otimes \rr/\lan a \otimes
(1-a):a, 1-a \in \rr \ran.
\end{equation}
See \cite[Corollary 4.3]{nes-suslin1990} and
\cite[Corollary 4.2]{mirzaii-2009} for different proofs. For the
isomorphism in the right hand side of (\ref{k2}) see also
\cite[Proposition 3.2.3]{guin1989}.

Important examples of rings with many units are semilocal rings which their
residue fields are infinite. In particular for an infinite  field $F$,
any commutative finite dimensional $F$-algebra is a semilocal ring and so
it is a ring with many units.

\begin{rem}
For a ring $R$ with many units and for any $n \ge 1$,
there exist $n$ elements in $R$ such that the sum of each nonempty
subfamily belongs to $\rr$ \cite[Proposition 1.3]{guin1989}.
Rings with this property are considered by Nesterenko
and Suslin \cite{nes-suslin1990}. For more about rings with many units
we refer to \cite{vdkallen1977}, \cite{nes-suslin1990}, \cite{guin1989} and
\cite{mirzaii2008}.
\end{rem}

The {\it pre-Bloch group} $\ppp(R)$ of a ring $R$ is defined as the quotient group
of the free abelian group $Q(R)$ generated by symbols $[a]$, $a, 1-a \in \rr$,
to the subgroup generated by elements of the form
\[
[a] -[b]+\bigg[\frac{b}{a}\bigg]
- \bigg[\frac{1- a^{-1}}{1- b^{-1}}\bigg]
+ \bigg[\frac{1-a}{1-b}\bigg],
\]
where $a,1-a, b, 1-b, a-b \in \rr$. Define
$\lambda': Q(R) \arr \rr \otimes \rr$, given by
$[a] \mt a \otimes (1-a)$. Then
\[
\lambda'\bigg([a] -[b]+\bigg[\frac{b}{a}\bigg]
- \bigg[\frac{1- a^{-1}}{1- b^{-1}}\bigg]
+ \bigg[\frac{1-a}{1-b}\bigg]\bigg)=
a \otimes \bigg( \frac{1-a}{1-b}\bigg)+\bigg(\frac{1-a}{1-b}\bigg)\otimes a.
\]
Let $(\rr \otimes \rr)_\sigma :=\rr \otimes \rr/
\lan a\otimes b + b\otimes a: a, b \in \rr \ran$.
We denote the elements of $\ppp(R)$ and $(\rr \otimes \rr)_\sigma$
represented by $[a]$ and $a\otimes b$ again by $[a]$ and $a\otimes b$,
respectively. Hence we have a well-defined map
\[
\lambda: \ppp(R) \arr (\rr \otimes \rr)_\sigma, \ \ \
[a] \mt a \otimes (1-a).
\]
The kernel of this map is called the {\it Bloch group} of $R$ and is
denoted by $B(R)$. Therefore for a ring  $R$ with many units,
using (\ref{k2}), we obtain the exact sequence
\begin{equation}
0 \arr B(R) \arr \ppp(R) \overset{\lambda}{\arr} (\rr \otimes \rr)_\sigma
\arr K_2^M(R) \arr 0.
\end{equation}

\section{Third homology of $\GL_2$}

Let $D_h(R^2)$ be the free $\z$-module with a basis consisting of
$(l+1)$-tuples $(v_0, \dots, v_l)$, where every
$\min\{l+1, 2\}$ of $v_i \in R^2$
are a basis of a direct summand of $R^2$.
Let us define a differential operator
\[
\partial_l : D_l(R^2) \arr D_{l-1}(R^2), \ \ l\ge 1,
\]
as an alternating sum of face operators $d_i$, which throw away the
$i$-th component of generators and let
$\partial_0: D_0(R^2) \arr\z$, $\sum_i n_i( v_i ) \mt \sum_i n_i$.
It is well known that the complex
\begin{gather*}
D_\ast \ : \ \ \
\cdots \larr  D_2(R^2) \overset{\partial_2}{\larr}
D_1(R^2) \overset{\partial_1}{\larr}
D_0(R^2) \overset{\partial_0}{\larr} \z \larr 0
\end{gather*}
is exact \cite[Lemma 1]{mirzaii2008}.

Let $H_{1}(Y):=\ker(\partial_{1})$. (See Remark 1.1 in
\cite{mirzaii2007} for an explanation for the choice of
the notation). From the short exact sequence
\[
0\larr \partial_{3}(D_{3}(R^2)){\larr} D_2(R^2){\larr} H_{1}(Y) \larr 0
\]
one obtains the exact sequence
\begin{gather*}
H_1(\GL_2,D_2(R^2))\arr T(R) \arr \PP(R)
\arr H_0(\GL_2,D_2(R^2)) \arr S(R) \arr 0,
\end{gather*}
where
\begin{gather*}
\hspace{-0.9 cm}
S(R):=H_0(\GL_2, H_1(Y)),\\
\PP(R):=H_0(\GL_2, \partial_3(D_{3}(R^2))), \\
\hspace{-0.9 cm}
T(R):=H_1(\GL_2, H_1(Y)).
\end{gather*}
By the Shapiro lemma \cite[Chap. III, Proposition~6.2]{brown1994},
$H_1(\GL_2,D_2(R^2))=0$ and
\[
H_0(\GL_2,D_2(R^2))=D_{2}(R^2)_{\GL_2}
\simeq \coprod_{a_i \in \rr} \z.[a_1, a_2].
\]
Therefore we get the exact sequence
\begin{equation}\label{bloch-exact1}
0 \larr T(R) \larr \PP(R) \overset{\nu}{\larr} \coprod_{a_i \in \rr}
\z.[a_1,a_2] \larr S(R) \larr 0.
\end{equation}

\begin{lem}
The group $\PP(R)$ is isomorphic to the quotient group of
the free abelian group generated by the symbols
$\left[
\begin{array}{ccc}
a_1 & a_2    \\
\lambda_1 & \lambda_2
\end{array}
\right]$,
$a_i$, $\lambda_i$, $1-\lambda_i$, $\lambda_1- \lambda_2 \in \rr$,
to the subgroup generated by the elements
\begin{gather*}
\hspace{-5 cm}
\left[
\begin{array}{cc}
(\lambda_2-\lambda_1)a_2 & \lambda_1   \\
\frac{\gamma_2-\gamma_1}{\lambda_2- \lambda_1}& \frac{\gamma_1}{\lambda_1}
\end{array}
\right]
-
\left[
\begin{array}{cc}
(\lambda_1-\lambda_2)a_1 & \lambda_2   \\
\frac{\gamma_1-\gamma_2}{\lambda_1- \lambda_2}& \frac{\gamma_2}{\lambda_2}
\end{array}
\right]
\\
\hspace{4.5 cm}
+
\left[
\begin{array}{cc}
\lambda_1a_1 & \lambda_2a_2   \\
\frac{\gamma_1}{\lambda_1}& \frac{\gamma_2}{\lambda_2}
\end{array}
\right]
-
\left[
\begin{array}{cc}
a_1 & a_2   \\
\gamma_1& \gamma_2
\end{array}
\right]
+
\left[
\begin{array}{cc}
a_1 & a_2   \\
\lambda_1& \lambda_2
\end{array}
\right],
\end{gather*}
where $\lambda_i, 1-\lambda_i, \gamma_i, 1-\gamma_i,
\lambda_1-\lambda_2, \gamma_1-\gamma_2, \lambda_i-\gamma_j \in \rr$.
\end{lem}
\begin{proof}
By applying the functor $H_0$ to
$D_{4}(R^2)\arr D_{3}(R^2) \arr \partial_{3}(D_{3}(R^2)) \arr 0$
we get the exact sequence
$D_{4}(R^2)_{\GL_2} \arr D_{3}(R^2)_{\GL_2} \arr \PP(R)\arr 0$.
It is easy to see that
\begin{gather*}
D_{3}(R^2)_{\GL_2}\simeq \coprod_{a, \lambda} \z.p(a,\lambda),\ \ \
D_{4}(R^2)_{\GL_2}\simeq \coprod_{a, \lambda,\gamma} \z.p(a,\lambda, \gamma),
\end{gather*}
where $p(a,\lambda)$ and $p(a,\lambda, \gamma)$ are the orbits of the frames
\begin{gather*}
\begin{array}{ll}
(e_1,  e_2 , a_1 e_1+a_2 e_2, \lambda_1 a_1 e_1+\lambda_2 a_2 e_2) &
\in D_{3}(R^2),\\
(e_1,  e_2 , a_1 e_1+a_2 e_2, \lambda_1 a_1 e_1+\lambda_2 a_2 e_2,
\gamma_1 a_1 e_1+\gamma_2 a_2 e_2) & \in D_{4}(R^2),
\end{array}
\end{gather*}
respectively, where $\lambda_i, 1-\lambda_i, \gamma_i, 1-\gamma_i,
\lambda_1-\lambda_2, \gamma_1-\gamma_2, \lambda_i-\gamma_j \in \rr$.
Set
\[
\left[
\begin{array}{cc}
a_1 &  a_2    \\
\lambda_1 & \lambda_2
\end{array}
\right]
:=p(a,\lambda) \mod D_{4}(R^2)_{\GL_2},
\]
which is an element of
$D_{3}(R^2)_{\GL_2}/D_{4}(R^2)_{\GL_2} \simeq \PP(R)$.
By a direct computation of $\partial_4(p(a,\lambda, \gamma))=0$
in $D_3(R^2)_{\GL_2}/D_4(R^2)_{\GL_2}$ one arrives at triviality of
the elements that is mentioned in the lemma.
\end{proof}

By a direct computation
\[
\nu(
\left[
\begin{array}{cc}
a_1 & a_2   \\
\lambda_1& \lambda_2
\end{array}
\right] )=[(\lambda_2-\lambda_1)a_2,
\lambda_1]-[(\lambda_1-\lambda_2)a_1, \lambda_2]+
[\lambda_1a_1, \lambda_2a_2]-[a_1, a_2],
\]
where $\nu$ is the map in the exact sequence (\ref{bloch-exact1}).
Thus $S(R)$ can be considered as the quotient
group of the free abelian group generated by the symbols $[a_1,a_2]$,
$a_1,a_2 \in \rr$ to the subgroup generated by
\[
[(\lambda_2-\lambda_1)a_2,
\lambda_1]-[(\lambda_1-\lambda_2)a_1, \lambda_2]+
[\lambda_1a_1, \lambda_2a_2]-[a_1, a_2],
\]
where $a_i, \lambda_i, 1-\lambda_i, \lambda_1-\lambda_2 \in \rr$.

\begin{sublem}\label{map1}
Let $\psi: \coprod_{a_i \in \rr} \z.[a_1,a_2] \arr (\rr
\otimes\rr)_\sigma$ be defined by $[a,b] \mt a \otimes b$. Then
\begin{gather*}
\psi\circ\nu(
\left[
\begin{array}{cc}
a_1 & a_2   \\
\lambda_1& \lambda_2
\end{array}
\right] )
=
(1-\lambda_1/\lambda_2) \otimes \lambda_1/\lambda_2 -
(-\lambda_2)\otimes \lambda_2\\
\hspace{3.6 cm}
=-(1-\lambda_2/\lambda_1) \otimes
\lambda_2/\lambda_1 + (-\lambda_1)\otimes \lambda_1.
\end{gather*}
\end{sublem}
\begin{proof}
This is obtained by a direct computation.
\end{proof}

This sublemma implies that the map
$\phi: S(R) \arr K_2^M(R)$ defined by  $[a, b] \arr \{a, b\}$
is well-defined.

\begin{lem}\label{sr-tr}
$S(R)\simeq \z \oplus K_2^M(R)$ and
$T(R)\simeq H_3(\GL_2)/H_3(\GL_1)$.
\end{lem}
\begin{proof}
From the exact sequence
$0 \arr H_1(Y) \arr D_1(R^2) \arr D_0(R^2) \arr \z \arr 0$,
one obtains a first quadrant spectral sequence
\[
E_{p, q}^1=\begin{cases}
H_q(\GL_2, D_p(R^2)) & \text{if $p=0, 1, 2$ }\\
H_q(\GL_2, H_1(Y)) & \text{if $p=3$ }\\
0 & \text{if $p \ge 4,$}\end{cases}
\]
which converges to zero.
Using the Shapiro lemma \cite[Chap. III, Proposition~6.2]{brown1994}
and a theorem of Suslin \cite[Theorem~1.9]{suslin1985},
\cite[2.2.2]{guin1989},
\[
E_{p, q}^1\simeq H_q(\GL_{2-p}), \ \ p=0,1,2.
\]
It is not difficult to see that $d_{1, q}^1=H_q(\inc)$, for $p=1,2$
\cite[Lemma 2.4]{nes-suslin1990}.
Thus the $E^1$-terms of the spectral sequence is as follows

\begin{gather*}
\begin{array}{llcccc}
\ast       &              &        &         &        &  \\
H_3(\GL_2) & H_3(\GL_1)   &  0     &         &        &  \\
H_2(\GL_2) & H_2(\GL_1)   &  0     &   \ast  &  0     &  \\
H_1(\GL_2) & H_1(\GL_1)   &  0     &   T(R)  &  0     &  \\
\z         &      \z      & \z     &   S(R)  &  0     &  0.
\end{array}
\end{gather*}
An easy analysis of this spectral sequence gives us the exact sequence
\[
0 \arr H_2(\GL_2)/ H_2(\GL_1) \arr S(R) \arr \z \arr 0
\]
and the isomorphism $T(R)\simeq H_3(\GL_2)/ H_3(\GL_1)$.
Now we need to prove that $K_2^M(R)\simeq H_2(\GL_2)/ H_2(\GL_1)$.
This is a well known fact. But here we give a rather simple proof.

First note that $\SK_1(R):=\SL(R)/E(R)=0$, where $E(R)$ is the elementary
subgroup of $\GL$. This follows from the
homology stability theorem $K_1(R)= H_1(\GL)\simeq H_1(\GL_1)\simeq \rr$
\cite[Theorem 1]{guin1989} and the fact that
$K_1(R) \simeq \rr \times \SK_1(R)$.

From the corresponding Lyndon-Hochschild-Serre spectral
sequence of the extension
$1 \arr \SL \arr  \GL \overset{\det}{\arr} \rr \arr 1$, using
the fact that $\SK_1(R)=0$, we obtain the decomposition
$H_2(\GL) \simeq H_2(\rr) \oplus K_2(R)$. To obtain this
decomposition we use the fact that
$E(R)$ is a perfect group, and  $K_2(R):=H_2(E(R))$.
Now by the homology stability theorem $H_2(\GL_2)\simeq H_2(\GL)$,
\cite[Theorem 1]{guin1989},
we have $K_2(R)\simeq H_2(\GL_2)/ H_2(\GL_1)$.
Thus by (\ref{k2}) we have
$K_2^M(R) \simeq H_2(\GL_2)/ H_2(\GL_1).$
Clearly, $\z \arr S(R)$ given by  $1 \mt [1,1]$
splits the above exact sequence. It is not difficult to see
that the projection $S(R) \arr K_2^M(R)$
is given by $\phi:[a,b] \mt \{a,b\}$.
This complete the proof of the lemma.
\end{proof}

\begin{rem}
It should be mentioned that
the isomorphism $S(R)\simeq \z \oplus K_2^M(R)$ in Lemma
\ref{sr-tr} already is proven in
\cite{nes-suslin1990} and \cite{guin1989}. But
their results are more general and so their methods are difficult.
\end{rem}

The map that gives the isomorphism $T(R)\simeq H_3(\GL_2)/ H_3(\GL_1)$
can be constructed directly. From the short exact sequence
\[
0 \arr \partial_1(D_1(R^2)) \arr D_0(R^2) \arr \z \arr 0
\]
we get the connecting homomorphism
$H_3(\GL_2) \arr H_2(\GL_2,\partial_1(D_1(R^2)))$.
Iterating this process
we get a homomorphism $\rho: H_3(\GL_2) \arr T(R)$.
Since the epimorphism $D_0(R^2) \arr \z$ has a
$\GL_{1}$-equivariant section $m \mt m( e_2)$, the
restriction of $\rho$ to $H_3(\GL_1)$ is zero. Thus we obtain
a homomorphism
\[
H_3(\GL_2)/H_3(\GL_1) \arr T(R).
\]
This map gives the mentioned isomorphism.

\begin{thm}\label{bloch-wigner-like}
There is an exact sequence
\begin{gather*}\label{pre-bloch-wigner}
0 \larr H_3(\GL_2)/H_3(\GL_1) \larr \PP(R)
\overset{\nu}{\larr}\!\!\!\!\!\!\!\!\!
\underset{(a_1,a_2) \neq (1,1)} {\coprod_{a_i \in \rr}}\!\!\!\!\!\!\!\!\!
\z.[a_1,a_2]
\larr K_2^M(R) \larr 0.
\end{gather*}
\end{thm}
\begin{proof}
As we mentioned in the proof of  Lemma
\ref{sr-tr} the map $\z \arr S(R)$, in the decomposition of
$S(R)$, is given by $1 \mt [1,1]$.
Now the exact sequence easily follows
from the exact sequence (\ref{bloch-exact1}) and Lemma \ref{sr-tr}.
We should mention that here
$\nu(\left[
\begin{array}{cc}
1 &  1   \\
\lambda_1 & \lambda_2
\end{array}
\right])=[\lambda_2-\lambda_1,
\lambda_1]-[\lambda_1-\lambda_2, \lambda_2]+
[\lambda_1, \lambda_2]$
and otherwise it is as above.
\end{proof}

\section{Connection to bloch-wigner exact sequence}

Consider the evident homomorphism
\[
\theta: \PP(R) \larr \ppp(R), \ \ \
\left[
\begin{array}{ccc}
a_1 & a_2    \\
\lambda_1 & \lambda_2
\end{array}
\right]
\mt [\lambda_1/\lambda_2].
\]
Sublemma \ref{map1} shows that the diagram
\[
\xymatrix{
\PP(R)    \ar[r]^{\!\!\!\!\!\!\!\!\!\!\!\!\!\!\!\!\!\!\!\!\!\! \nu}
\ar[d]^{\theta} &
\coprod_{a_i \in \rr} \z.[a_1,a_2] \ar[d]^\psi \\
\ppp(R) \ar[r]^{\!\!\!\!\!\!\!\!\!\!\!\!\!\!\!\!\!\!\!\!\! \lambda}
& (\rr \otimes\rr)_\sigma}
\]
dose not commute in general. But it commutes if we restrict $\nu$ to its kernel.
Thus $\theta$ maps $H_3(\GL_2)/H_3(\GL_1)$
into $B(R)$. It is well known \cite{suslin1991}, \cite{mirzaii-2009} that
this map, i.e. $H_3(\GL_2)/H_3(\GL_1)\arr B(R)$, is surjective, which relies on
some tedious calculation. We denote this map again by $\theta$.
The natural inclusion $\rr \times \GL_1 \arr \GL_2$, gives us the
cup product map
\[
\cup: \rr \otimes H_2(\GL_1) \arr H_3(\GL_2), \ \
a \otimes(b \wedge c) \mt a \cup (b \wedge c).
\]
It is not difficult to see that the composition
\[
\rr \otimes H_2(\GL_1) \overset{\cup}{\arr} H_3(\GL_2) \arr
H_3(\GL_2)/H_3(\GL_1) \overset{\theta}{\arr} B(R)
\]
is trivial. So we obtain the map $\H_3(\SL_2(R), \z) \arr B(R)$,
where
\[
\H_3(\SL_2(R), \z):= H_3(\GL_2)/\im\Big(H_3(\GL_1)+ \rr \cup H_2(\GL_1)\Big).
\]
The following theorem is the main theorem of \cite{mirzaii-2009}.

\begin{thm}\label{bloch-wigner1}
For a commutative ring $R$ with many units, we have the exact sequence
\[
\tors(\mu_R, \mu_R) \arr \H_3(\SL_2(R), \z) \arr B(R) \arr 0.
\]
When $R$ is an integral domain, then the left hand side map in the above
exact sequence is injective.
\end{thm}
\begin{proof}
See Theorem 5.1 of \cite{mirzaii-2009}.
\end{proof}

Motivated by Sublemma \ref{map1}, we define the map
\[
\eta: \ppp(R) \larr (\rr \otimes \rr)_{\overline{\sigma}}, \ \ \
[a] \mt a \otimes (1-a),
\]
where $(\rr \otimes\rr)_{\overline{\sigma}}:=
\rr \otimes\rr/\lan  a\otimes (-a): a \in \rr \ran$. Since
elements of the form $a\otimes b + b \otimes a$ vanish in
$(\rr \otimes\rr)_{\overline{\sigma}}$, $\eta$ is well-defined.
We denote $\ker(\eta)$ by $B'(R)$.
Thus we have the exact sequence
\[
0 \larr B'(R) \larr \ppp(R) \overset{\eta}{\larr}
(\rr \otimes \rr)_{\overline{\sigma}} \larr K_2^M(R) \larr 0.
\]

Let
\[
\psi':\coprod_{a_i \in \rr} \z.[a_1,a_2] \larr
(\rr \otimes \rr)_{\overline{\sigma}}
\]
be defined by $[a_1,a_2] \arr a_1\otimes a_2$ and let
$\alpha: (\rr \otimes \rr)_{{\sigma}} \larr
(\rr \otimes \rr)_{\overline{\sigma}}$
be the canonical map.
Note that we have the following commutative diagram with exact rows
\[
\xymatrix{
0 \ar[r] & T(R)
\ar[r] \ar[d]& \PP(R)
\ar[r]^{\!\!\!\!\!\!\!\!\!\!\!\!\!\!\!\!\!\!\!\!\!\!\!\!\!\!\nu} \ar[d]^\theta &
\underset{(a_1,a_2)\neq (1,1)}{\coprod_{a_i \in \rr}} \z.[a_1,a_2]
\ar[r] \ar[d]^{\psi'}& K_2^M(R)\ar[r] \ar[d]^{=} & 0 \\
0 \ar[r] & B'(R) \ar[r] &
\ppp(R) \ar[r]^\eta & (\rr \otimes\rr)_{\overline{\sigma}} \ar[r] &
K_2^M(R)\ar[r]  & 0 \\
0 \ar[r] & B(R) \ar[r] \ar[u] &
\ppp(R) \ar[r]^\lambda  \ar[u]_{=} &
(\rr \otimes\rr)_{{\sigma}} \ar[r]\ar[u]_\alpha&
K_2^M(R)\ar[r] \ar[u]_{=} & 0.
}
\]

In the rest of this section we study the group $B'(R)/B(R)$.
Clearly  $B(R)$ is a subgroup of $B'(R)$.
Let $\lan a \ran:=[a]+[a^{-1}] \in \ppp(R)$, which is a $2$-torsion element
\cite[Lemma 1.2]{suslin1991}.
By a direct computation,
\[
\lambda\Big(\lan a \ran\Big)=a \otimes(-a) \in (\rr \otimes\rr)_{{\sigma}}
\]
and
\[
\lambda\Big(\lan ab \ran-\lan a \ran -\lan b \ran\Big)=a \otimes
b + b \otimes a= 0 \in (\rr \otimes\rr)_{{\sigma}}.
\]
Thus $\lan a \ran \in B'(R)$ and
$\lan ab \ran-\lan a \ran -\lan b \ran \in B(R)$.

\begin{prp}\label{prp}
There is a surjective map
\[
\rr/\lan -1, r^2: r \in \rr \ran \larr B'(R)/B(R).
\]
Moreover $B'(R)/B(R)$ is a  $2$-torsion group generated by the elements
$\lan a \ran$.
\end{prp}
\begin{proof}
To prove this we look at the following commutative diagram
\[
\xymatrix{
0 \ar[r] & B(R)
\ar[r] \ar[d]& \ppp(R)
\ar[r]^{\!\!\!\!\!\!\!\!\!\!\!\!\!\!\!\lambda}
\ar[d]^{=} &
(\rr \otimes\rr)_{\sigma}
\ar[r] \ar[d]^{\alpha}& K_2^M(R)\ar[r] \ar[d]^{=} & 0 \\
0 \ar[r] & B'(R) \ar[r] &
\ppp(R) \ar[r]^{\!\!\!\!\!\!\!\!\!\!\!\!\!\!\!\eta}
& (\rr \otimes\rr)_{\overline{\sigma}} \ar[r] &
K_2^M(R)\ar[r]  & 0.}
\]
Breaking this diagram into two diagrams with short exact sequence rows
and applying the Snake lemma, one obtains the isomorphism
\[
\ker(\alpha) \simeq B'(R)/B(R).
\]
Clearly $\ker(\alpha)$, as a subgroup of
$(\rr \otimes\rr)_{\overline{\sigma}}$, is generated by the elements
$a \otimes (-a)$.
Consider the map $\rr \arr \ker(\alpha)$, given by
$a \mt a \otimes (-a)$. This map is surjective and
$\lan -1, r^2: r \in \rr \ran$ is in its kernel.

Let $a \in \rr$. Consider the surjective linear forms defined on $R^2$
by the vectors $f_1=(1,0)$, $f_2=(0,1)$, $f_3=(1,-1)$, $f_4=(1,-a)$. Since
$R$ is a ring with many units,  there exists $(x, y) \in R^2$ such that
$f_i(x,y) \in \rr$, $i=1, \dots, 4$. Therefore there exist $c \in R$ such that
$c, 1-c, 1-ac \in \rr$. Since
\[
a \otimes (-a)=ac \otimes (-ac) - c \otimes (-c)- a \otimes c - c \otimes a,
\]
the group $\ker(\alpha)$ is generated be the elements $a \otimes (-a)$, where
$a, 1-a \in \rr$. Now by a direct computation one can see that
under the isomorphism $\ker(\alpha) \simeq B'(R)/B(R)$, the element
$a \otimes (-a)$, $a, 1-a \in \rr$, maps to $\lan a \ran$.
(See also Example \ref{example}(iii) below).
\end{proof}

\begin{exa}\label{example}
(i) If $\rr=\rr^2$, then $B(R)=B'(R)$. For example of rings with this property 
see \cite[Corollary 5.6]{mirzaii-2009}.
\par (ii) Clearly $\R^\times=\lan -1, r^2: r \in \R^\times \ran$. Therefore
$B(\R)=B'(\R)$.
\par (iii) Let $F$ be an infinite field and let
$\tau: \fff \arr \ppp(F)$ and $\varsigma: \fff \arr (\fff \otimes\fff)_{{\sigma}}$
are given by $a \mt \lan a \ran$ and $a \mt a \otimes(-a)$,
respectively.
Then the following diagram commutes
\[
\xymatrix{
\fff/\fff^2 \ar[r]^{\!\!\!\!\!\!\!\!\!=} \ar [d]^\tau
& \fff/\fff^2 \ar [d]^\varsigma\\
\ppp(F) \ar[r]^{\!\!\!\!\!\!\!\!\!\!\!\!\lambda} & (\fff \otimes\fff)_{\sigma}.}
\]
Thus the surjective map in Proposition \ref{prp}
\[
\fff/\lan -1, r^2: r \in \fff \ran \larr B'(F)/B(F)
\]
is given by $a \arr \lan a \ran$.
For any ring with many units the map $\tau$ may not be defined!
\end{exa}

\subsection*{Acknowledgements}
Part of this work was done during my stay at IH\'ES. I would like
to thank them for their support and hospitality.


\bigskip
\address{{\footnotesize

Department of Mathematics,

Institute for Advanced Studies in Basic Sciences,

P. O. Box. 45195-1159, Zanjan, Iran.

email:\ bmirzaii@iasbs.ac.ir
}}

\end{document}